\documentclass[12pt]{amsart}

\usepackage{amsmath}
\usepackage{amsfonts}
\usepackage{amssymb}
\usepackage{latexsym}
\usepackage[dvips]{graphicx}
\newtheorem{theorem}{Theorem}[section]
\newtheorem{proposition}[theorem]{Proposition}
\newtheorem{corollary}[theorem]{Corollary}
\newtheorem{lemma}[theorem]{Lemma}
\newtheorem{definition}[theorem]{Definition}
\newtheorem{remark}[theorem]{Remark}
\theoremstyle{definition}

\numberwithin{equation}{section}

   \def\sH{{\mathfrak H}}

      \def\dC{{\mathbb C}}
\def\dD{{\mathbb D}}

   \def\dN{{\mathbb N}}   
      \def\dR{{\mathbb R}}

\def\cD{{\mathcal D}}      
   \def\cH{{\mathcal H}}

\def\cS{{\mathcal S}}

\def\RE{{\rm Re\,}}
\def\IM{{\rm Im\,}}
\def\ker{{\rm ker\,}}

\def\f{\varphi}

\def\ovl{\overline}

\def\ran{{\rm ran\,}}
\def\dom{{\rm dom\,}}

\usepackage{color}

\newcommand{\rd}{\textcolor{red}}

\begin{document}
{\qquad \hfill \rd{\textit{\textbf{{Draft-fin (Oct.16,2008)}}}}}

\vskip 3.5cm

\title[Numerical Range and Quasi-Sectorial Contractions]
{Numerical Range and Quasi-Sectorial Contractions}
\author{
Yury Arlinski\u{i}}
\address{Department of Mathematical Analysis \\
East Ukrainian National University \\
Kvartal Molodyozhny 20-A \\
Lugansk 91034 \\
Ukraine} \email{yury$\_$arlinskii@yahoo.com; yma@snu.edu.ua}
\author{Valentin Zagrebnov}
\address{Universit\'{e} de la M\'editerran\'ee  and  Centre de Physique
Th\'{e}orique - UMR 6207 \\  Luminy - Case 907, Marseille 13288, Cedex 9, France}
\email{zagrebnov@cpt.univ-mrs.fr}



\thispagestyle{empty}

\begin{abstract}
We apply a method developed by one of the authors, see
\cite{Arl1}, to localize the numerical range of \textit{quasi-sectorial} contractions semigroups.
Our main theorem establishes a relation between the numerical range of quasi-sectorial contraction
semigroups $\{\exp(- t S)\}_{t\ge 0}$, and the maximal {sectorial} generators $S$.
We also give a new prove of the rate $O(1/n)$ for the operator-norm Euler formula approximation:
$\exp(- t S)=\lim\limits_{n\to \infty}(I+tS/n)^{-n}$, $t\ge 0$, for this class of semigroups.
\end{abstract}
\maketitle

\tableofcontents

\section{Introduction}

In what follows the Banach algebra of  all \textit{bounded} linear
operators on a complex Hilbert space $\sH$ is denoted by
$\mathcal{L}(\sH)$. We denote by $I_\sH$ the identity operator in a
Hilbert space $\sH$. The \textit{domain}, the \textit{range}, and
the \textit{null-space} of a linear operator $T$ are denoted by
$\dom T$, $\ran T,$ and $\ker T$, respectively. For
$T\in\mathcal{L}(\sH)$ the operators $\RE T = (T + T^*)/2$ and $\IM
T = (T - T^*)/2i$ are said to be the \textit{real} and the
\textit{imaginary} parts of $T$.

\subsection{Numerical range and the Kato mapping theorem}

Let $\sH$ be a complex separable Hilbert space and let $A$ be an
(\textit{unbounded}) linear operator in $\sH$ with domain $\dom A$.
\begin{definition}\label{NUMRANGE}
The set of complex numbers
\[
W(A):=\left\{(Au,u) \in \mathbb{C}: \; u\in\dom A,\; ||u||=1\right\}
\]
is called the \textit{numerical range} of $A$, or its field of values.
\end{definition}
According to the Hausdorff-Toeplitz theorem, the numerical range is
a convex set. We recall also the following properties of the numerical
range (see e.g. \cite{Ka}).
\begin{proposition}\label{Prp}
Let $A$ be a closed operator in $\sH$. Then
\begin{enumerate}\def\labelenumi{\rm (\alph{enumi})}
\item
for any complex number $z\notin\overline{W(A)}$ holds $\ker
(A-zI_\sH)=\{0\}$ and $\ran (A-zI_\sH)$ is closed. Moreover, the
defect
\[
{\rm def}(A-zI_{\sH}):=\dim\left(\sH\ominus\ran (A-zI_\sH)\right)
\]
is constant in each connected component of
$\dC\setminus\overline{W(A)}$.
\item
If $z\in\dC\setminus\overline{W(A)}$ then
\[
||(A-zI_{\sH})^{-1} f||\le\cfrac{1}{{\rm dist}(z,W(A))} \ ||f||,\;
f\in\ran (A-zI_\sH).
\]
\item
If $\dom A$ is dense in $\sH$ and $W(A)\ne \dC$, then $A$ is
closable.
\end{enumerate}
\end{proposition}
{\begin{corollary}\label{spect-n-ran-bounded}
For a bounded operator $A\in \mathcal{L}(\mathfrak{H})$ the spectrum $\sigma(A)$ is a subset of
$\ \overline{W(A)}$.
\end{corollary}
For unbounded operator $A$ the relation between spectrum and numerical range is more complicated.
We would like to warn, that it may very well happen that $\sigma(A)$ is not contained
in $\overline{W(A)}$, but for a closed operator $A$ the essential spectrum
$\sigma_{ess}(A)$ is always a subset of $\overline{W(A)}$. The condition
${\mathrm{def}}(A-zI) = 0, \ z\notin \overline{W(A)}$ in Proposition \ref{Prp}
serves to ensure that for those unbounded operators one gets
\begin{equation}\label{spect-n-ran}
\sigma(T)\subset \overline{W(A)} \ ,
\end{equation}
i.e., the same conclusion as in Corollary \ref{spect-n-ran-bounded} for bounded operators.}

In the sequel we need the following numerical range mapping theorem
due to  Kato \cite{Kato}.
\begin{proposition} \label{nummap} \cite{Kato}.
Let $f(z)$ be a rational function on $\dC$, with $f(\infty)=\infty$.
Let for some compact and convex set $E'\subset\dC$ the inverse
function $f^{-1}:E'\to E \, \supseteq \, K$, where $K$ is a convex
kernel of $E$, i.e., is a subset of $E$ such that $E$ is star-shaped
with respect to any $z\in K$.
If $A$ is bounded operator with $W(A)\subseteq K$, then $W(f(A)) \subseteq E'$.
\end{proposition}
\noindent Notice that for a convex set $E$ the \textit{maximal} convex kernel $K = E$.


\subsection{Sectorial operators and quasi-sectorial contractions}

\begin{definition}
\label{accret}\cite{Ka}. Let $S$ be a linear operator in a Hilbert
space $\sH$. If $\RE(Su,u)\ge 0$ for all $u\in\dom S$, then $S$ is
called accretive.
\end{definition}
So, the operator $S$ is accretive if and only if its numerical range
is contained in the closed right-half plane of the complex plane. An
accretive operator $S$ is called \textit{maximal} accretive ($m$-accretive)
if one of the equivalent conditions is satisfied:
\begin{itemize}
\item the operator $S$ has no accretive extensions in $\sH$;
\item
the resolvent set $\rho(S)$ is nonempty;
\item the operator $S$ is densely defined and closed,  and $S^*$ is
accretive operator.
\end{itemize}
The resolvent set $\rho(S)$ of $m$-accretive operator contains the
open left half plane and
\[
||(S-zI_\sH)^{-1}||\le\cfrac{1}{|\RE z|} \ , \ \RE z<0.
\]
It is well known \cite{Ka} that if $S$ is $m$-accretive operator, then the
one-parameter semigroup
$$T(t)=\exp(-tS), \;t\ge 0$$
is contractive. Conversely, if the family $\{T(t)\}_{t\geq 0}$ is a strongly continuous
semigroup of bounded operators in a Hilbert space $\sH$, with
$T(0)=I_\sH$ ($C_0$-semigroup) and $T(t)$ is a contraction for each
$t$, then the generator $S$ of $T(t)$:
\[
Su:=\lim\limits_{t\to +0}\cfrac{(I_\sH-T(t)) u}{t}\ ,\ u\in\dom S \ , \]
where domain is defined by condition:
\[ \dom S=\left\{u\in\sH:\lim\limits_{t\to
+0}\cfrac{(I_\sH-T(t))u}{t}\quad\mbox{exists}\right\} \ ,
\]
is an $m$-accretive operator in $\sH$. Then the Euler formula approximation:
\begin{equation} \label{Euler}
T(t)=s-\lim\limits_{n\to\infty}\left(I_\sH+\cfrac{t}{n}\,
S\right)^{-n},\; t\ge 0
\end{equation}
holds in the strong operator topology, see e.g. \cite{Ka}.
\begin{definition}
\label{sect}\cite{Ka}. Let $\alpha\in [0,\pi/2)$ and let
\[
\mathcal{S}(\alpha):=\left\{z\in\dC:|\arg z|\le \alpha\right\}
\]
be a sector on the complex plane $\dC$ with the vertex at the origin
and the semi-angle $\alpha$.

A linear operator $S$ in a Hilbert space $\sH$ is called sectorial
with vertex at $z=0$ and the semi-angle $\alpha$ if $W(S)\subseteq
\mathcal{S}(\alpha)$.
\end{definition}
If $S$ is $m$-accretive and sectorial with vertex at $z=0$ and the
semi-angle $\alpha$ then it is called $m$-sectorial with
vertex at the origin and with the semi-angle $\alpha$. For short we call these operators
$m$-$\alpha$-sectorial. The resolvent set of $m$-$\alpha$-sectorial operator $S$ contains the set
$\dC\setminus \mathcal{S}(\alpha)$ and
\[
||(S-zI_\sH)^{-1}||\le\cfrac{1}{{\rm
dist}\left(z,\mathcal{S}(\alpha)\right)} \ ,\ z\in \dC\setminus \mathcal{S}(\alpha).
\]

It is well-known \cite{Ka} that
\begin{itemize}
\item
 a $C_0$-semigroup $T(t)=\exp(-tS),$ $t\ge 0$ has
\textit{contractive} and \textit{holomorphic} continuation into the
sector $\mathcal{S}(\pi/2-\alpha)$ if and only if the generator $S$
is $m$-$\alpha$-sectorial operator,
\item the \textit{sesquilinear} form: $(Su,v) \ ,\ u,v\in\dom S$ is
closable.
\end{itemize}
Denote by $\mathfrak{S}[\cdot,\cdot]$ the closure of the form $(Su,v)$ and by
$\cD[\mathfrak{S}]$ its domain. By the \textit{first representation
theorem} \cite{Ka} the operator $S$ is \textit{uniquely} associated with the form
$\mathfrak{S}[\cdot,\cdot]$ in the following sense:
\[
(Su, v)=\mathfrak{S}[u,v]\quad\mbox{for all}\quad u\in\dom S\quad\mbox{and for
all}\quad v\in\cD[\mathfrak{S}].
\]

The important relation between the one-parameter semigroups $\{T(t):=\exp(-tS)\}_{t\geq0}$
generated by the $m$-$\alpha$-sectorial operators $S$ and the corresponding
closed sesquilinear forms is established in \cite{Arl3}:
\[
u\in\cD[\mathfrak{S}]= \left\{u\in \mathfrak{H}:
\frac{d}{dt}(T(t)u,u)\Bigl|_{t=+0} \ \mbox{exits}\right\},
\]
and one has that
\begin{equation}\label{semi-gr-deriv}
\frac{d}{dt}(T(t)u,u)\Bigl|_{t=+0}=-\mathfrak{S}[u,u] \ ,\ u\in\cD[\mathfrak{S}] \ .
\end{equation}

Notice that in this case of  the Euler approximation \eqref{Euler}
converges to the semigroup in the \textit{operator-norm} topology \cite{Cach},
\cite{CaZag}, \cite{Zagrb2}.
\begin{definition}\label{alpha-domain} {\cite{CaZag}}.
For any $\alpha\in[0, \pi/2)$ we define in the complex plane
${\mathbb{C}}$ a closed domain:
\begin{equation}\label{domain-D}
D_{\alpha}:=\{z\in {\mathbb{C}}: |z|\leq \sin \alpha\} \cup
\{z\in {\mathbb{C}}: |\arg (1-z)|\leq \alpha \ {\rm{and}}\ |z-1|\leq
\cos \alpha \} \ .
\end{equation}
This is a convex subset of the unit disc $\dD = D_{\alpha =\pi/2}$, with the "angle"
{\rm{(}}in contrast to the "tangent"{\rm{)}} touching of the boundary $\partial \dD$
at the only one point $z=1$. It is evident that $D_{\alpha}\subset D_{\beta > \alpha}$.
\end{definition}
\begin{definition}\label{Q-S-cont}{\cite{CaZag}}.
A contraction $C$ on the Hilbert space $\mathfrak{H}$ is called
\textit{quasi-sectorial} with  semi-angle $\alpha\in [0, \pi/2)$, if its numerical range
$W(C) \subseteq D_{\alpha}$.
\end{definition}

It is evident that if operator $C$ is a \textit{quasi-sectorial} contraction, then
$I - C$ is an $m$-\textit{sectorial} operator with
vertex $z=0$ and semi-angle $\alpha$. The limits $\alpha=0$ and $\alpha = \pi/2$ correspond,
respectively, to non-negative (i.e. \textit{self-adjoint}) and to some \textit{general} contraction.
\begin{remark}\label{res-fam}\cite{CaZag}.
{Notice that the resolvent family $\{(I_\sH + t \, S)^{-1}\}_{t\geq
0}$ of the $m$-$\alpha$-sectorial operator $S$,  gives the first
non-trivial example of a quasi-sectorial contractions, if one
considers the semi-angles $\alpha\in [0, {\pi/3})$. Below (see
Section 2) we show that it can be extended to $\alpha\in [0,
\pi/2)$.}
\end{remark}
The definition of quasi-sectorial contractions was motivated in
\cite{CaZag} by the lifting of the Trotter-Kato product formula and
the Chernoff theory of semigroup approximation \cite{Chern} to the
operator-norm topology.

Namely, using properties of quasi-sectorial contractions established
in \cite{CaZag} for the case of $m$-$\alpha$-sectorial generator, it
was proved that there is \textit{operator-norm} convergence of the
Euler semigroup approximation \eqref{Euler} and of the Trotter
product formula, see \cite{Cach}, \cite{CaZag}, \cite{P}, \cite{BP},
\cite{Zagrb2}. Theorem 2.1 from \cite{CaZag} states that the
operators $T(t)=\exp(-tS)$ are quasi-sectorial contractions with
$W(T(t))\subset D_\alpha$ for all $t\ge 0 $ and $\alpha\in [0,
\pi/2)$. Here operator $S$ stands for $m$-$\alpha$-sectorial
generator. As it is indicated in \cite{Zagrb1}, the proof in
\cite{CaZag} has a flaw. In \cite{Zagrb1} it was corrected, but only
for the semi-angles $\alpha\in [0,{\pi}/{4}]$.

Our main Theorem \ref{main} establishes a quite accurate relation
between $m$-$\alpha$-sectorial generators and the numerical range of
the corresponding one-parameter contraction semigroups. It improves
the recent result \cite{Zagrb1} from  $\alpha\in [0, \pi/4]$ to
$\alpha\in [0, \pi/2)$. To this end we use in the next section some
nontrivial results due to \cite{Arl1}-\cite{Arl4}, concerning a class
of operator contractions and of semigroups on the complex plane.
Besides that we use recent results related to the generalizations of the
famous von Neumann inequality \cite{vN}, obtained in \cite{BC},
\cite{C}, \cite{CD}, \cite{DD}, in order to give a new prove that in fact the Euler formula \eqref{Euler}
converges in the operator-norm with the rate $O(1/n)$.

\section{Operators of the class $C_\sH(\alpha)$}

{A fundamental for us will be the class of contractions introduced for the first time in
\cite{Arl1} and studied in \cite{Arl1}-\cite{Arl2}:}
\begin{definition}\label{CA} \cite{Arl1}.
Let $\alpha\in(0,\pi/2)$. We say that a bounded operator $T \in \mathcal{L}(\mathfrak{H})$
belongs to the class $C_\sH(\alpha)$ if
\begin{equation}
\label{DEFCA} ||T \, \sin\alpha  \pm i I_\sH \cos\alpha ||\le 1.
\end{equation}
\end{definition}
It is clear, that the class $C_\sH(\alpha)$ is a convex and closed
(with respect to the strong operator topology) set, which is
intersection of two closed operator balls corresponding to $\pm$.
Moreover, by virtue of (\ref{DEFCA}) one immediately concludes that:
\[
T\in C_\sH(\alpha)\iff -T\in C_\sH(\alpha)\iff T^*\in C_\sH(\alpha),
\]
and that condition \eqref{DEFCA} is equivalent to the following criterium:
\begin{equation}\label{DEFCA-bis}
\tan\alpha\, (||f||^2-||Tf||^2)\ge 2|\IM(Tf,f)|,\; f\in \sH.
\end{equation}
{This inequality implies that the operator $T$ is a contraction.
Together with Definition \ref{NUMRANGE}, it also proves that $T\in
C_\sH(\alpha)$ is equivalent to the statement that $(I-T^*)(I+T)$ is
a bounded $m$-$\alpha$-sectorial operator.}

According to (\ref{DEFCA-bis}) it is natural to identify
$C_\sH(\alpha=0)$ with the set of all \textit{self-adjoint}
contractions whereas $C_\sH(\alpha=\pi/2)$ is the set of general
contractions on $\sH$.

Now following \cite{Arl4} we define for $\alpha\in(0,\pi/2)$ the family of subsets of the complex plane $\mathbb{C}$:
\begin{equation}\label{CCA}
C(\alpha)=\left\{z\in \dC:|z\sin\alpha\pm i\cos\alpha|\le
1\right\}=\left\{z\in\dC:  (1-|z|^2)\tan\alpha \geq 2|\IM z| \right\},
\end{equation}
cf. definitions (\ref{DEFCA}), (\ref{DEFCA-bis}). Then, similar to (\ref{DEFCA}), each set $C(\alpha)$ is the intersection
of \textit{two} closed disks of the complex plane which is contained in the closed unit disk
$\ovl\dD=\{z\in\dC:|z|\le 1\}$. From (\ref{CCA}) it is clear that $C(\alpha=0)=[-1,1].$
\begin{remark}\label{NR-mapping}
By virtue of definitions  (\ref{DEFCA}), (\ref{DEFCA-bis}) and
(\ref{CCA}) one gets that $T\in C_\sH(\alpha)$ implies for the
numerical range: $W(T)\subseteq C(\alpha)$.  Then one has:
\begin{equation}
\label{IMPART} ||\IM T||\le\tan ({\alpha}/{2}) \ .
\end{equation}
Notice that $C(\alpha)\setminus D_{\alpha}\neq \emptyset$,
i.e. in general the operator $T$ is not quasi-sectorial.
\end{remark}
Besides, the operator class $C_\sH(\alpha)$ has several interesting
properties. In particular the following one \cite{Arl1}:
\begin{proposition} \label{Comm}\cite{Arl1}. Let $T_1$ and $T_2$
belong to the class $C_\sH(\alpha)$. Then the operator
\[
(T_1T_2+T_2T_1)/2
\]
 belongs to the class $C_\sH(\alpha)$. In
particular, if $T_1$ and $T_2$ are two commuting operators from the
class $C_\sH(\alpha)$, then the product $T = T_1 \, T_2$ also
belongs to $C_\sH(\alpha)$.
\end{proposition}
\begin{remark}\label{multiplic-semi}
Notice that the set $C(\alpha)$ inherits a similar property: for each $\alpha\in(0,\pi/2)$ it forms a multiplicative
semigroup of complex numbers in the plane $\dC$ . Detailed properties of these semigroups have been studied
in \cite{Arl4}.
\end{remark}

{Now we are in position to establish a connection between set of contractions $C_\sH(\alpha)$
and $m$-$\alpha$-sectorial operators.}
\begin{proposition}\label{class-C-sectorial}\cite{Arl1}.
If $S$ is $m$-$\alpha$-sectorial operator, then
\[
T=(I_\sH-S)(I_\sH+S)^{-1}
\]
 belongs to the class $C_\sH(\alpha)$ and
conversely, if $T\in C_\sH(\alpha)$, and $\ker(I_\sH+T)=\{0\}$, then
$S=(I_\sH-T)(I_\sH+T)^{-1}$ is $m$-$\alpha$-sectorial operator.
\end{proposition}
Now let $S$ be $m$-$\alpha$-sectorial operator and let $\lambda>0$.
Then $\lambda S$ is also $m$-$\alpha$-sectorial. By Proposition
\ref{class-C-sectorial} the operator $U(\lambda)=(I_\sH-\lambda
S)(I_\sH+\lambda S)^{-1}$ belongs to the class $C_\sH(\alpha)$. Let
us put
\begin{equation}\label{F-T}
F(\lambda):= (I_\sH+\lambda S)^{-1}=\frac{1}{2}(U(\lambda)+I_\sH).
\end{equation}
Since the operators $U(\lambda)$ and $I_\sH$ belong to the class
$C_\sH(\alpha)$ and since it is a convex set, the operator
$F(\lambda)$ is also in the class $C_\sH(\alpha)$. Hence, by Remark
\ref{NR-mapping} one obtains $W(F(\lambda))\subseteq C(\alpha)$.
\begin{remark}\label{L-alpha}
Following the arguments of \cite{CaZag},\cite{Zagrb1} we can
localize the numerical range $W(F(\lambda))$ even more accurate, cf.
Remark \ref{res-fam}. Since for any $u\in\mathfrak{H}$ we have:
\begin{equation*}
(u, F(\lambda) u) = (v_\lambda, v_\lambda) + \lambda (S v_\lambda, v_\lambda) \in \mathcal{S}(\alpha) \ ,
\end{equation*}
where $v_\lambda:= F(\lambda) u$, it follows that
\begin{equation*}
W(F(\lambda))\subseteq\mathcal{S}(\alpha) \cap C(\alpha) \ ,
\end{equation*}
for $\lambda> 0$ and $\alpha\in [0,\pi/2)$. Moreover, since $U(\lambda)=2F(\lambda)-I_\sH \in
C_\sH(\alpha)$, see (\ref{F-T}), we find that
\begin{equation}
\label{NUMC}
\left\|\left(F(\lambda) - I_\sH /2 \right) \sin\alpha \pm i I_\sH
(\cos\alpha)/2 \right\|\le {1}/{2} \ .
\end{equation}
For the numerical range this implies:
\begin{equation}\label{ident-2}
W(F(\lambda))\subseteq L(\alpha):= \{\zeta\in \dC: |(\zeta -
1/2)\sin\alpha \pm (i\cos\alpha)/2 |\le 1/2\}\ ,
\quad\mbox{for}\quad \lambda> 0 \ .
\end{equation}
Notice that $L(\alpha)\subset\mathcal{S}(\alpha) \cap C(\alpha)$, see Figure 1.
\end{remark}

{Now, we follow essentially the line of reasoning of \cite{Arl1} to establish the one-to-one correspondence
between $m$-$\alpha$-sectorial generators and contraction semigroups of the class  $C_\sH(\alpha)$, cf.
Remark \ref{NR-mapping} about relation to quasi-sectorial contractions.}
\begin{theorem} \label{SEMCA}
If $S$ is $m$-$\alpha$-sectorial operator in a Hilbert space $\sH$,
then the corresponding semigroup $T(t)=\exp(-tS)\in C_\sH(\alpha)$,
for all $t\ge 0$. Conversely, let $T(t)=\exp(-tS)$ for $t\ge 0$ be a
$C_0$-semigroup of contractions on a Hilbert space $\sH$. If
$T(t)\in C_\sH(\alpha)$ for non-negatives $t$ in neighborhood of $t=0$,
then the generator $S$ is an $m$-$\alpha$-sectorial operator.
\end{theorem}
\begin{proof} Let $S$ be a $m$-$\alpha$-sectorial operator and let $\lambda\ge 0$.
By (\ref{F-T}) the
operator $F(\lambda)=(I_\sH+\lambda S)^{-1}$ belongs to the class
$C_\sH(\alpha)$. Therefore, by Proposition \ref{Comm} for each $t\ge 0$ and any natural
number $n$ the operator
\[
F^n\left(\frac{t}{n}\right)=\left(I_\sH+\frac{t}{n}S\right)^{-n}
\]
belongs to the class $C_\sH(\alpha)$. Taking in account that the set
$C_\sH(\alpha)$ is closed with respect to the strong operator
topology, from the Euler formula \eqref{Euler} we get that strong
limit $T(t)=\exp(-tS)$ also belongs to the class $C_\sH(\alpha)$.

Now suppose that semigroup $T(t)=\exp(-tS)\in C_\sH(\alpha)$ for  $t\in [0,\delta )$,
where $\delta>0$. Define operator family:
\[
B_{\pm}(t):= T(t) \, \sin\alpha\pm i\cos\alpha I_\sH, \; t\ge 0.
\]
Since $B_\pm(0)=(\sin\alpha\pm i\cos\alpha) I_\sH$ and $T(t)\in
C_\sH(\alpha)$ for $t\in [0,\delta)$, we get
\[
||B_\pm (t)f||^2\le ||f||^2=||B_\pm(0)||f||^2,\; t\in [0,\delta),\;
f\in \sH.
\]
Since
\[
||B_\pm(t)f||^2=\sin^2\alpha ||T(t)f||^2+\cos^2\alpha ||f||^2\pm
2\sin\alpha\,\cos\alpha\,\IM (T(t)f,f),\; f\in\sH,
\]
for all $f\in\dom S$ we have:
\[
\frac{d}{dt}\left(||B_\pm(t)f||^2\right)\Bigl|_{t=+0}=-2\sin^2\alpha\,
\RE(Sf,f)\pm 2\sin\alpha\,\cos\alpha\,\IM (Sf,f)\le 0.
\]
Thus, $W(S)\subseteq \mathcal{S}(\alpha)$. But since operator $S$ is $m$-accretive, it is
$m$-$\alpha$-sectorial \cite{Ka}.
\end{proof}

\section{Numerical range for contractive holomorphic semigroups}

{From Theorem \ref{SEMCA} it follows, in particular, that for
$m$-$\alpha$-sectorial generator $S$ the numerical range of the
corresponding contraction semigroup
\[
W(\exp(-tS))\subseteq C(\alpha)\quad\mbox{for all}\quad t\ge 0 \ .
\]
But as we warranted in Remark \ref{NR-mapping} it does not imply that this semigroup is quasi-sectorial
contraction.
It was discovered in \cite{CaZag} that the conformal mapping : $z \mapsto z^2$ , together with the Kato
numerical range
theorem (Proposition \ref{nummap}) play a special r\^{o}le in the theory of quasi-sectorial contractions.}
\begin{definition} \label{Omega}
Let $\alpha\in [0,\pi/2)$. We define a domain:
\begin{equation}\label{Omega1}
\Omega(\alpha):=\left\{z^2: z\in C(\alpha)\right\}.
\end{equation}
\end{definition}
So, if $f(z)=z^2$, then $\Omega(\alpha)=f(C(\alpha)).$ Since (see Remark \ref{multiplic-semi}) $C(\alpha)$ is a
\textit{multiplicative semigroup}, we obtain that $\Omega(\alpha)\subseteq C(\alpha)$, and that the
subset $\Omega(\alpha)$ is in turn a multiplicative semigroup \cite{Arl4}.

{From (\ref{CCA}) and Proposition \ref{Comm} it
follows then that for any $\alpha\in(0,\pi/2)$ the set (\ref{Omega}) has representation:
\begin{equation}\label{OM}
\Omega(\alpha)=\left\{z\in\dC: |\sqrt{z}\sin\alpha\pm
i\cos\alpha|\le
1\right\}=\left\{z\in\dC:2|\IM\sqrt{z}|\le(1-|z|)\tan\alpha \right\} ,
\end{equation}
with the limiting cases: $\Omega(\alpha = 0)=[0,1]$ and $\Omega(\alpha = \pi/2)= \ovl\dD$.}
\begin{lemma} \label{conv}
{The set $\Omega(\alpha)$ {\rm{(}}\ref{Omega1}{\rm{)}} is convex and
$\Omega(\alpha)\subseteq D_{\alpha}$.}
\end{lemma}
\begin{proof} Let $C_{+}(\alpha):=\{z\in C(\alpha):  \IM z\ge 0\}$. Then
clearly, $\Omega(\alpha)=f(C_{+}(\alpha))$, where $f(z)=z^2$. Denote
\[
\Gamma(\alpha)=\partial C_{+}(\alpha)\setminus(-1,1).
\]
Then $\partial\Omega(\alpha)=f(\Gamma(\alpha))$. Since
\[
\Gamma(\alpha)=\left\{z:z=\frac{e^{it}-i\cos\alpha}{\sin\alpha},\;
t\in\left [\frac{\pi}{2}-\alpha,\frac{\pi}{2}+\alpha\right]\right\},
\]
the boundary $\partial\Omega(\alpha)$ can be parameterized as
follows:
\begin{equation}\label{boundary-Omega}
\partial\Omega(\alpha)=\left\{z=\zeta(t)=\frac{(e^{it}-i\cos\alpha)^2}{\sin^2\alpha},\;
t\in\left [\frac{\pi}{2}-\alpha,\frac{\pi}{2}+\alpha\right]\right\}.
\end{equation}
Put $x=x(t):=\RE\zeta(t)$, $y=y(t):=\IM\zeta(t)$. Then, since
\[
\frac{d^2 y}{dx^2}=\cfrac{y''(t)x'(t)-y'(t)x''(t)}{(x'(t))^3} \ ,
\]
we obtain
\[
\frac{d^2
y}{dx^2}=-\cfrac{\IM(\zeta'(t)\ovl{\zeta''(t))}}{(\RE\zeta'(t))^3} \ .
\]
Further, by explicit calculations we get:
\[
\zeta'(t)=\frac{2ie^{it}(e^{it}-i\cos\alpha)}{\sin^2\alpha} \ \ , \ \
\zeta''(t)=\frac{2ie^{it}(2ie^{it}+\cos\alpha)}{\sin^2\alpha} \ ,
\]
\[
- \IM(\zeta'(t)\ovl{\zeta''(t)})=
\cfrac{4(2+\cos^2\alpha-3\cos\alpha\sin
t)}{\sin^4\alpha} \ \ , \ \
\RE\zeta'(t)=\cfrac{2\cos t(\cos\alpha-2\sin t)}{\sin^2\alpha} \ \, .
\]
For $t\in[{\pi}/{2}-\alpha , {\pi}/{2}+\alpha]$, by the estimates
\[
\cos^2\alpha-3\cos\alpha\sin t+2\ge 2 \cos^2-3\cos\alpha
+2=(2-\cos\alpha)(1-\cos\alpha)>0,
\]
we obtain $-\IM(\zeta'(t)\ovl{\zeta''(t)})>0.$

Since  $\cos\alpha-2\sin t\le -\cos\alpha<0$, then for
$t\in[{\pi}/{2}-\alpha , {\pi}/{2}+\alpha]$, we get that
\[
\begin{array}{l}
\cfrac{d^2 y}{dx^2}<0 \quad\mbox{for}\quad t\in\left
[\cfrac{\pi}{2}-\alpha , \cfrac{\pi}{2}\right) \ ,\\
\cfrac{d^2 y}{dx^2}>0 \quad\mbox{for}\quad t\in
\left(\cfrac{\pi}{2} , \cfrac{\pi}{2}+\alpha\right] \ ,
\end{array}
\]
which implies that the set $\Omega(\alpha)$ is convex.

Since the mapping: $z \mapsto z^2$ is conformal and $\Omega(\alpha)$ is convex, the proof of the second part
of the lemma follows from the  estimate
\begin{equation}\label{est-D}
- (\tan \alpha/2)^2 \leq \RE\zeta(t) \ ,
\end{equation}
where $(\tan \alpha/2)^2 \leq \sin \alpha $ for
$\alpha\in(0,\pi/2)$,  see (\ref{domain-D}) and
(\ref{boundary-Omega}).
\end{proof}

Notice that in view of relation:
\[
y(t)=\frac{\sin 2t-2\cos\alpha\cos t}{\sin^2\alpha} \ ,
\]
one obtains
\begin{equation}\label{maxmod}
\max\limits_{z\in\Omega(\alpha)}|\IM z|=\frac{\sin
2\gamma-2\cos\alpha\sin\gamma}{\sin^2\alpha} \ ,
\end{equation}
where
\begin{equation}
\label{extrm}
\sin\gamma=\frac{\cos\alpha+\sqrt{\cos^2\alpha+8}}{4} \ , \ \;\gamma\in
(0,\pi/2) \ .
\end{equation}
Since $\Omega(\alpha)\subset C(\alpha)$,  the number in the
right-hand side of \eqref{maxmod} is less than $\tan
({\alpha}/{2}).$
\begin{lemma}\label{LM1}
Let $\alpha\in [0,\pi/2)$ and let $T\in C_\sH(\alpha)$. Then $W(T^{2n})\subseteq\Omega(\alpha)$ for all
natural numbers $n$.
\end{lemma}
\begin{proof} We apply the Kato numerical range mapping theorem (Proposition
\ref{nummap}) for $f(z)=z^2$ with $E'=\Omega(\alpha)$. Then by Definition \ref{Omega} we have
$E:=f^{-1}(\Omega(\alpha))=C(\alpha)$. Since $C(\alpha)$ is a convex
set, its maximal convex kernel $K$ \textit{coincides} with
$C(\alpha)$. By virtue of Proposition \ref{Comm} from $T\in C_\sH(\alpha)$ it
follows that for all natural $n$ one gets $T^n\in C_\sH(\alpha)$. Then by Remark
\ref{NR-mapping} we obtain $W(T^n)\subseteq C(\alpha) = K$. Applying now the Kato
mapping theorem for $f(z)=z^2$ we obtain $W(T^{2n})\subseteq\Omega(\alpha)$.
\end{proof}

Now we are in position to proof the main theorem of the present paper.
{\begin{theorem}\label{main}
{\rm{(1)}} Let $S$ be $m$-$\alpha$-sectorial operator. Then
\begin{equation}\label{num-ran-main}
W(\exp(-tS))\subseteq \Omega(\alpha) \ , \ t\ge 0.
\end{equation}
In particular, $\{\exp(-tS)\}_{t\geq0}$ is the quasi-sectorial
contraction semigroup.

{\rm{(2)}} The inverse is also true. Let $\{T(t):=\exp(-tS)\}_{t\geq0}$ be a
$C_0$-semigroup on a Hilbert space $\sH$. If in some neighborhood of
$t=0$ the numerical range: $\,W(T(t\geq0))\subseteq\Omega(\alpha)$ for some $\alpha\in[0,\pi/2)$,
then the generator $S$ is an $m$-$\alpha$-sectorial operator.
\end{theorem}
\begin{proof}
(1) By Theorem \ref{SEMCA} the operators $T(t)=\exp(-tS)$ belong to
the class $C_\sH(\alpha)$ for all $t\ge 0$, but $T(t)=T^2(t/2)$, then
by Lemma \ref{LM1} we obtain $W(T(t))\subseteq\Omega(\alpha)$. Since Lemma \ref{conv} implies
$\Omega(\alpha)\subset D_\alpha$, the operators $T(t)$ are quasi-sectorial contractions.

(2) Fix $u\in \dom S$, $||u||=1$. For $t\in [0,\delta]$ we define two functions:
\[
G_\pm(t):=|\sin\alpha\,\sqrt{(T(t)u,u)}\pm i\cos\alpha|^2 \ .
\]
Then $G_\pm(0)=1$. The condition $W(T(t))\subseteq\Omega(\alpha)$
yields (see \eqref{OM}) that $G_\pm (t)\le G_\pm(0)$, $t\in [0,\delta]$, which implies:
\[
\frac{d}{dt}\left(G_\pm(t)\right)\Bigl|_{t=+0}\le 0 \ .
\]
Therefore we get that
\begin{equation}\label{etim-below}
\lim\limits_{t\to+0}\frac{1-G_\pm (t)}{t}\ge 0.
\end{equation}
Further we use the following identities:
\[
\begin{array}{l}
1-G_\pm (t)=\sin^2\alpha (1-|(T(t)u,u)|)\mp
2\sin\alpha\,\cos\alpha\,\IM\sqrt{(T(t)u,u)}=\\
=\cfrac{\sin^2\alpha
\left(1-(T(t)u,u)+(T(t)u,u)(1-\overline{(T(t)u,u)})\right)}{1+|(T(t)u,u|}\mp
2\sin\alpha\,\cos\alpha\,\IM\sqrt{(T(t)u,u)},
\end{array}
\]
\[
\begin{array}{l}
\cfrac{1-G_\pm (t)}{t}=\cfrac{\sin^2\alpha
\left(1-(T(t)u,u)+(T(t)u,u)(1-\overline{(T(t)u,u)})\right)}{t(1+|(T(t)u,u|)}\mp\\
\mp 2
\sin\alpha\,\cos\alpha\,\IM\left(\cfrac{(T(t)u,u)-1}{t(\sqrt{(T(t)u,u)}+1)}\right) \ .
\end{array}
\]
Since by (\ref{semi-gr-deriv}) one has:
\[
\lim\limits_{t\to+0}\frac{1-(T(t)u,u)}{t}=(Su,u)\ ,
\]
and by (\ref{etim-below}) we get
\[
\lim\limits_{t\to+0}\frac{1-G_\pm (t)}{t}=\sin^2\alpha \,
\RE(Su,u)\pm\sin\alpha\,\cos\alpha\,\IM(Su,u)\ge 0 \ .
\]
The last estimate implies that $W(S)\subseteq\cS(\alpha)$, i.e., $S$ is
$m$-$\alpha$-sectorial operator.
\end{proof}
\noindent Thus, the two equivalent conditions: $\exp(-tS)\in
C_\sH(\alpha)$ (see Theorem \ref{SEMCA}) and
$W(\exp(-tS))\subseteq\Omega(\alpha)$ \textit{completely} characterize
$m$-$\alpha$-sectorial operators $S$. Since (\ref{num-ran-main})
yields $W(\exp(-tS))\subset D_\alpha$, the statement in Theorem
\ref{main}(1) is extension of Theorem 2.1 (proven in \cite{Zagrb1}
for $\alpha\in [0,\pi/4)$) to the whole class of the
$m$-$\alpha$-sectorial generators: $\alpha\in [0,\pi/2)$.

Notice that taking into account the multiplicative semigroup property of the set $\Omega(\alpha)$ and
Theorem \ref{main} we get by consequence the following inequalities:
\[
\left|\sin\alpha\,\left(\prod\limits_{k=1}^n\left(\exp(-t_k
S_k)u_k,u_k)\right)\right)^{1/2}\pm i\cos\alpha\right|\le 1
\]
for arbitrary $m$-$\alpha$-sectorial operators $S_1,\ldots S_n$,
any non-negative numbers $ t_1,\ldots, t_n,$ and arbitrary
normalized vectors $u_1,\ldots,u_n$ from $\sH$. By semi-group
property: \\
$\exp(-t_1 S)\exp(-t_2 S)=\exp(-(t_1+t_2)S)$. Hence we also
get
\[
|\sin\alpha\sqrt{(\exp(-t_1 S)u, \exp(-t_2 S^*)u)}\pm
i\cos\alpha|\le 1
\]
for any $t_1, t_2\ge 0$ and any $u\in\sH$, $||u||=1$.}

Remark that since $\RE z\ge -\tan^2({\alpha}/{2})$ for all
$z\in\Omega(\alpha)$ and $W(\exp(-tS))\subseteq\Omega(\alpha)$, we obtain:
\[
 \RE \exp(-tS)\ge-\tan^2({\alpha}/{2}) \ I_\sH \ , \ t\ge 0 \ .
\]
The inclusion $W(\exp(-tS))\subseteq\Omega(\alpha)$ also implies
(see \eqref{maxmod}) that
\[
||\IM \exp(-tS)||\le\frac{\sin
2\gamma-2\cos\alpha\sin\gamma}{\sin^2\alpha} < \tan({\alpha}/{2}) \ , \  t\ge 0 \ ,
\]
where $\gamma$ satisfies \eqref{extrm}. By virtue of relation $\RE T^2=(\RE T)^2-(\IM T)^2$
(valid for any bounded operator $T$) and since one has that
$\RE \exp(-tS)=(\RE\exp(-tS/2))^2-(\IM\exp(-tS/2))^2$, we obtain:
\[
\RE \exp(-tS)\ge - \left(\frac{\sin
2\gamma-2\cos\alpha\sin\gamma}{\sin^2\alpha}\right)^2
I_\sH \ > - \tan^2({\alpha}/{2}) \ I_\sH \ , \ t\ge 0 \ .
\]

Now it is useful to define also the following subset of $C(\alpha)$:
\begin{equation}
\label{Q}
\begin{array}{l}
Q(\alpha):=\left\{z\in C(\alpha):|z\sin\alpha-\cos\alpha|\le
1\right\}=\\
=\{z\in\dC:2|\IM z|\le (1-|z|^2)\tan\alpha\}\cap \{z\in\dC:-2\RE z\le(1-|z|^2)\tan\alpha\} \ .
\end{array}
\end{equation}
So, $Q(\alpha)$ is the intersection of \textit{three} closed disks, cf. (\ref{CCA}).
\begin{proposition}
\label{SemQ} The set $Q(\alpha)$ is a multiplicative semigroup on $\mathbb{C}$.
\end{proposition}
\begin{proof} 
{First we note that $Q(\alpha)$ contains the set
$C_{+}(\alpha):= \{z\in C(\alpha): \RE z\ge 0$\}.
Moreover, since $\RE iz = - \IM z$, by (\ref{Q}) we obtain also that the set
\[
B(\alpha):=\left\{z\in C(\alpha):iz\in C(\alpha)\right\}  \ ,
\]
has a non-empty intersection with $Q(\alpha)$, such that}
\begin{equation}\label{prop}
B(\alpha) \cap Q(\alpha) = B(\alpha) \ \  \textrm{and} \ \
Q(\alpha)\setminus B(\alpha) \subset C_{+}(\alpha) \ .
\end{equation}
{In \cite{Arl4} it is shown  that the set $B(\alpha)$ is the \textit{ideal} of
the multiplicative semigroup $C(\alpha)$, i.e.,  $z\, \xi \in B(\alpha)$
for all $z\in B(\alpha)$ and all $\xi\in C(\alpha)$. Let $z,\xi\in
Q(\alpha)\setminus B(\alpha)$. Since $B(\alpha)$ has the properties (\ref{prop}), then}
\[
\arg z,\arg\xi\in (- {\pi}/{4}, {\pi}/{4} ).
\]
Finally, since $z\,\xi\in C(\alpha)$ and $\arg(z\,\xi)\in
(- {\pi}/{2}, {\pi}/{2} )$, we obtain $z\,\xi\in Q(\alpha).$
\end{proof}

Notice that in particular Proposition \ref{SemQ} yields:
$$z\in
Q(\alpha)\Rightarrow z^n\in Q(\alpha).$$ Moreover, if $\f_n(z)=z^n$,
then the set $\f_n(Q(\alpha))$ is also a multiplicative semigroup and
$$\f_n(Q(\alpha))\subset Q(\alpha).$$
Since by (\ref{CCA}) and (\ref{Q}) one gets that
\[
C(\alpha)=Q(\alpha)\cup\{-Q(\alpha)\},
\]
where $\{-Q(\alpha)\}:=\{z: -z\in Q(\alpha)\}$, by (\ref{Omega}) and (\ref{prop}) we have
\[
\Omega(\alpha)=f(Q(\alpha))
\]
when $f(z)=z^2$.  Thus, we obtain:
\begin{equation} \label{INCL} \Omega(\alpha)\subset
Q(\alpha).
\end{equation}
\begin{remark}\label{another way}
{Another way to check \eqref{INCL} is the following argument. Let
$\xi=z^2$, where $z\in C(\alpha)$. Then $\xi\in \Omega(\alpha)$, see
(\ref{Omega}). Since $\RE\xi=(\RE z)^2-(\IM z)^2$ and $(\IM
z)^2\le|\IM z|$, we also get}
\[
\begin{array}{l}
-2\RE\xi=-2(\RE z)^2+2(\IM z)^2\le 2|\IM z|\le
\tan\alpha(1-|z|^2)=\\
\qquad=\tan\alpha(1-|\xi|)\le\tan\alpha(1-|\xi|^2) \ ,
\end{array}
\]
which means by (\ref{Q}) that $\Omega(\alpha)\subset Q(\alpha)$.
\end{remark}
\begin{remark}\label{domains}
Let $D_\alpha$ be the set (\ref{domain-D}) introduced in
\cite{CaZag}. From definition (\ref{Q}) one deduces that
$Q(\alpha)\subset D_\alpha$, as well as that $D_\alpha$ is not a
subset of $C(\alpha)$, i.e., $D_\alpha\cap C(\alpha)\ne D_\alpha$.
In addition, by virtue of (\ref{ident-2}) and (\ref{OM}) we obtain:
$L(\alpha)\subset\Omega(\alpha)$, where $L(\alpha)$ is defined by
(\ref{ident-2}). So, for all $\alpha\in [0,\pi/2)$ and $F(\lambda)=
(I_\sH+\lambda S)^{-1},$ $\lambda>0$ we have the following
inclusions {\rm{(}}see Figure 1{\rm{):}}
\[
\begin{array}{l}
L(\alpha)\subset\Omega(\alpha)\subset Q(\alpha)\subset D_\alpha \ ,\\
W(F(\lambda))\subseteq L(\alpha) \ ,\\
W(F^n(\lambda))\subseteq C(\alpha)\quad\mbox{and}\quad
W(F^{2n}(\lambda))\subseteq \Omega(\alpha)\quad\mbox{for all}\quad n\in\dN \ ,
\end{array}
\]
\end{remark}

\includegraphics[width=\textwidth,keepaspectratio]{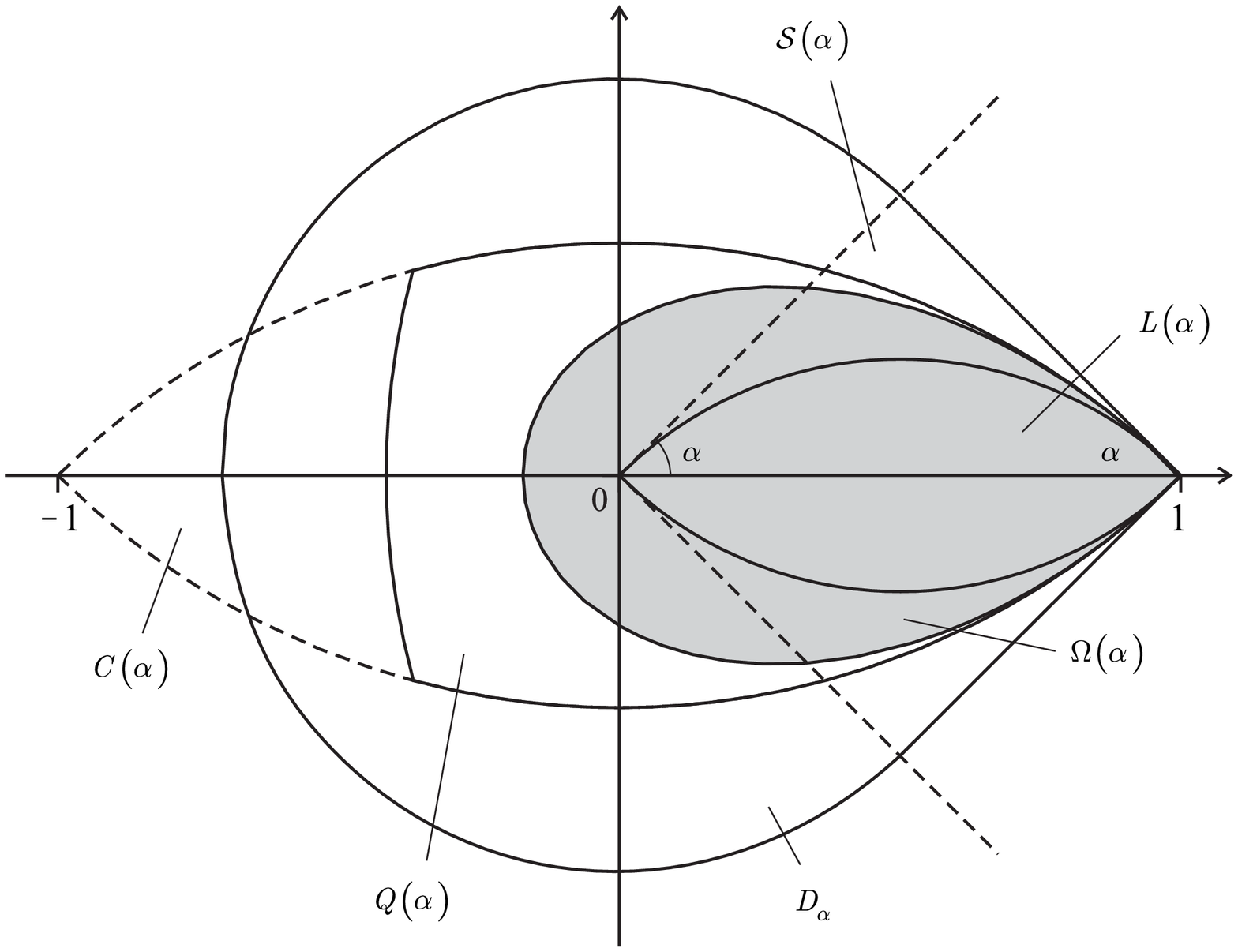}
\centerline{Figure 1}

\section{{Operator-norm convergence of the Euler formula}}
\subsection{{A brief history}}

Theorem 5.1 in \cite{CaZag} states that
\begin{equation}
\label{CZag1}
\left\|\left(I_\sH + {t}S/{n}\right)^{-n}-\exp(-tS)\right\|=O\left({\ln n}/{n}\right)
\end{equation}
for $m$-$\alpha$-sectorial boundedly invertible operator $S$ in $\sH$
and for all $t\in S_{\pi/2-\alpha}$. The proof of \eqref{CZag1} is
essentially based on two estimates for a quasi-sectorial contraction
$C$, i.e. $W(C)\subseteq D_\alpha$, also established in \cite{CaZag}. The first is
\begin{equation} \label{CZag2}
\left\|C^n-C^{n+1}\right\|\le \frac{K}{n+1}
\end{equation}
for all $n\in\dN$ and for some $K>0$ depending on $\alpha$, and the second one is
\begin{equation}\label{CZag3}
\left\|C^n-\exp{(n(C-I_\sH))}\right\|=O\left({1}/{n^{1/3}}\right) \ .
\end{equation}
Then application of \eqref{CZag2} and \eqref{CZag3} for the particular case of the operator
$C=(I_\sH + {t}S/{n})^{-1}$ leads to \eqref{CZag1}.

Later the estimate \eqref{CZag1} was improved in \cite{P}  to
\begin{equation} \label{P1}
\left\|C^n-\exp{(n(C-I_\sH))}\right\|=O\left(\left({\ln n}/{n}\right)^{1/2}\right),
\end{equation}
and in \cite{BP} to
\begin{equation} \label{BP1}
\left\|C^n-\exp{(n(C-I_\sH))}\right\|=O\left({1}/{n^{1/2}}\right).
\end{equation}
In \cite{P}  the following bound is proved for the operator-norm convergence rate:
\begin{equation}
\label{P2}
\left\|\left(I_\sH + {t}S/{n}\right)^{-n}-\exp(-tS)\right\| \le {c}/{n},
\end{equation}
where $c$ is a constant depending on the operator $S$. The bounds
\eqref{CZag3}, \eqref{P1}, \eqref{BP1}, and \eqref{P2} are
obtained in \cite{P} and \cite{BP} by means of the probability theory methods
(the \textit{Poisson distribution} and the \textit{Central Limit Theorem}) that improves the estimate
(\ref{CZag3}).

Here, using the results obtained in \cite{DD}, \cite{BC}, \cite{CD},
and \cite{C}, and  generalizations of some von
Neumann inequality \cite{vN}, we prove the theorem, which makes more explicit the right-hand side of the
estimate (\ref{P2}).
\begin{theorem}
\label{NormEst} Let $S$ be an $m$-$\alpha$-sectorial operator in a
Hilbert space $\sH$. Then
\[
\left\|\left(I_\sH + {t}S/n \right)^{-n}-\exp(-tS)\right\|\le {K(\alpha)}/{(n\cos^2\alpha)} \ , \ t\ge 0 \ ,
\]
where
\[
({\pi}\sin\alpha)/{2\alpha}\le
K(\alpha)\le\min\left(2+{2}/{\sqrt{3}} \ , \ {(\pi-\alpha)}/{\alpha}\right) \ .
\]
\end{theorem}
\subsection{{The von Neumann inequality and its generalizations}}

The spectral sets theory was introduced by von Neumann \cite{vN}
in order to extend the functional calculus to the case of non-normal operators in
Hilbert spaces.
\begin{definition} \cite{vN}.
\label {SS} A set $\sigma\subset\dC$ is a spectral set of the
operator $A$ in a Hilbert space $\sH$ if it is closed and if for any
bounded rational function $u(z)$ on $\sigma$ one has
\[
||u(A)||\le \sup\limits_{z\in\sigma}|u(z)|.
\]
\end{definition}
\begin{proposition} \cite{vN}.
\label{VonNeum} A necessary and sufficient condition for one of the
domains
\[
|z-a|\le r,\; |z-a|\ge r,\; \RE(az)\ge b
\]
to be a spectral set of $A$ in $\sH$ is that
\[
||A-a I_\sH||\le r,\;||(A- aI_\sH)^{-1}||\le \frac{1}{r},\;
\RE(aT)\ge b.
\]
\end{proposition}
\begin{definition}
\label{Kspect} Let $\cD$ be an open convex subset of the complex
plane {\rm{(}}$\cD\ne\emptyset,\; \cD\ne\dC${\rm{)}} and let $A$ be a linear
operator in a Hilbert space $\cH$ with $W(A)\subset\overline{\cD}$.
The set $\cD$ is called \textit{K-spectral} set for the operator $A$, if
\[
||u(A)||\le K\sup\limits_{z\in\cD}|u(z)|
\]
for all rational functions $u(z)$ without pole in the spectrum of
$A$.
\end{definition}
Let $\cD_k:=\{|z-a_k|< r_k\}$, $k=1,2$ be two disks such that
$\cD_k\cap\cD_2=\{\xi_1,\xi_2\}$. Let
\[
\mathfrak{L}:=\cD_1\cap\cD_2.
\]
The set $\mathfrak{L}$ is said to be a \textit{convex lens-shaped} domain
\cite{BC}. Denote by $2\alpha\in(0,\pi)$ the angle of the lens $L$
at the vertices. The operator $A\in \mathcal{L}(\sH)$ is called \textit{of
the lenticular $\mathfrak{L}$-type} \cite{BC} if
\[
||A-a_1 I_\sH||\le r_1\quad\mbox{and}\quad||A-a_2 I_\sH||\le r_2.
\]
The next Proposition is established in \cite{BC}.
\begin{proposition}
\label{BCTh} Let $\mathfrak{L}$ be a convex lens-shaped domain of the complex
plane with angle $2\alpha$. There exists a best positive constant
$K(\alpha)$ such that the inequality
\begin{equation}
\label{BCIn} ||p(A)||\le K(\alpha)\sup\limits_{z\in \mathfrak{L}}|p(z)|,
\end{equation}
holds for all polynomials $p(z)$, for all operators $A\in \mathcal{L}(\sH)$
of $\mathfrak{L}$-type, and for all Hilbert spaces $\sH$. The constant
$K(\alpha)$ depends only on $\alpha$ and satisfies the inequality
\[
\frac{\pi}{2\alpha}\sin\alpha\le
K(\alpha)\le\min\left(2+\frac{2}{\sqrt{3}},\,\,\frac{\pi-\alpha}{\alpha}\right).
\]
\end{proposition}
Notice that by the Mergelyan Theorem \cite{Rudin}, the inequality
\eqref{BCIn} remains valid if $p$ is holomorphic in $\mathfrak{L}$
and continuous in $\overline{\mathfrak{L}}$.

\subsection{{Proof of the Theorem \ref{NormEst}}}
\begin{proof}
Let
\[
\begin{array}{l}
\cD_1=\left\{z\in\dC:\left|z-\cfrac{1}{2}-\cfrac{i}{2}\cot\alpha\right|<\cfrac{1}{2\sin\alpha}\right\},\\
\cD_{{2}}=\left\{z\in\dC:\left|z-\cfrac{1}{2}+\cfrac{i}{2}\cot\alpha\right|<\cfrac{1}{2\sin\alpha}\right\}.
\end{array}
\]
Then (see \eqref{ident-2}) $\mathfrak{L} = L(\alpha) \setminus \partial L(\alpha)$, i.e. interior of the set $L(\alpha)$.

Fix $t\ge 0$ and let
\[
C=F(t):=\left(I_\sH+tS\right)^{-1} \ .
\]
Then from (\ref{NUMC}) it follows that the operator $C$ is of $\mathfrak{L}$-type.

Since
\[
F\left({t}/{n}\right) = F(t) \{(n-1)F(t)/{n} + I_\sH /{n} \}^{-1} \ , \ n\in\dN \ ,
\]
we get also that
\[
F^n\left({t}/{n}\right)=C^n \{(n-1) C /{n} + I_\sH /{n} \}^{-n} \ .
\]
Put
\[
h_n(z):=\exp(1-1/z)-z^n\{(n-1) z /{n} + {1}/{n}\}^{-n} \ , \
z\in\mathfrak{L} \ .
\]
Because $\RE z > 0$ for all $z\in\mathfrak{L}$, the function
$h_n(z)$ is holomorphic in $\mathfrak{L}$ and continuous in
$\overline{\mathfrak{L}}=L(\alpha)$. Moreover, since
\[
\exp(-tS)=\exp(I_\sH-C^{-1}) \ ,
\]
we obtain that
\[
\exp(-tS)-\left(I_\sH + {t} S /{n}\right)^{-n}= h_n(C) \ .
\]

The fractional-linear conformal transformation $w \mapsto
{z=1}/(1+w)$ maps the sector $\cS(\alpha)$ onto $L(\alpha)$. So, let
\[
g_n(w):=h_n((1+w)^{-1})=\exp(-w)-\left(1 + {w}/{n}\right)^{-n} \ , \
w\in (\cS(\alpha)\setminus \partial\cS(\alpha)) \ .
\]
Then clearly:
\[
\sup\limits_{z\in\mathfrak{L}}|h_n(z)|=\sup\limits_{w\in(\cS(\alpha)\setminus \partial\cS(\alpha))}|g_n(w)|
=\sup\limits_{w\in\partial\cS(\alpha)}|g_n(w)| \ .
\]
Since
\[
\partial\cS(\alpha)=\{x\exp(-i\alpha),\;x\in\dR_+\}\cup
\{x\exp(i\alpha),\;x\in\dR_+\} \ ,
\]
we have to estimate the value of
\[
\sup\limits_{x\in\dR_+}|g_n(x\exp(\pm i\alpha)| \ .
\]

To this end we use the representation:
\[
g_n(x\exp(i\alpha))=-\int\limits_{0}^x
\frac{d}{ds}\left((1+ {se^{i\alpha}}/{n})^{-n} \ {\exp(-(x-s)e^{i\alpha})}\right)ds \ ,
\]
where
\[
\frac{d}{ds}\left((1+ {se^{i\alpha}}/{n})^{-n} \ {\exp(-(x-s)e^{i\alpha})}\right) =
\cfrac{se^{i\alpha}}{n} \ (1+ {se^{i\alpha}}/{n})^{-(n+1)} \ {\exp(-(x-s)e^{i\alpha})}  \ .
\]
By elementary inequalities:
\[
\begin{array}{l}
\left|1 + {se^{i\alpha}}/{n}\right|^{n+1}=\left(1+ {s^2}/{n^2}+ {2s\cos\alpha}/{n}\right)^{(n+1)/2} \ge
\left(1+ {(s\cos\alpha)}/{n}\right)^{n+1}\ge\\
\quad\le 1+ {((n+1)} s \cos\alpha)/n \ge 1+s\cos\alpha,
\end{array}
\]
we obtain
\[
{s} \ {\left|1+ {(s\cos\alpha)}/{n}\right|^{-(n+1)}}\le {s}/{(1+s\cos\alpha)}\le
{1}/{\cos\alpha} \ .
\]
Therefore, we obtain as an upper bound:
\[
|g_n(x\exp(i\alpha))|\le
\int\limits_{0}^x\frac{\exp(-(x-s)\cos\alpha)}{n \cos\alpha}
ds\le\frac{1}{n\cos^2\alpha}(1-\exp(-x\cos\alpha))
\le\frac{1}{n\cos^2\alpha}.
\]
Similarly one obtains the estimate:
\[
|g_n(x\exp(-i\alpha))|\le \frac{1}{n\cos^2\alpha}.
\]
Thus,
\[
\sup\limits_{x\in\dR_+}|g_n(x\exp(\pm i\alpha)| \le
\frac{1}{n\cos^2\alpha}.
\]
Now by Proposition\ref{BCTh} we get
\begin{equation}\label{alpha-estim}
\left\|\exp(-t S)-\left(I_\sH+ {t} S/n\right)^{-n}
\right\|=||h_n(C)||\le \frac{K(\alpha)} {n \cos^2\alpha} \ ,\ t\ge 0 \ ,
\end{equation}
which completes the proof.
\end{proof}


\section{Conclusion}

\noindent Now several remarks are in order:

(a) Theorem 2.1 from \cite{CaZag} states that for quasi-sectorial
contractions one has:
\begin{equation}
\label{CAHZAG} W(\exp(-tS))\subseteq D_\alpha,\; \alpha \in
[0,{\pi}/{2}),
\end{equation}
Here $S$ stands for $m$-$\alpha$-sectorial generator of contraction.
As it is indicated in \cite{Zagrb1}, the proof in \cite{CaZag} has a
flaw. In \cite{Zagrb1} it was corrected but only for the range
$\alpha\in [0,{\pi}/{4}]$. Because
\[
\Omega(\alpha)\subset Q(\alpha)\subset D_\alpha,
\]
our main Theorem \ref{main} shows that the original claim
\eqref{CAHZAG} in \cite{CaZag} is indeed correct for all $\alpha\in
[0,{\pi}/{2})$.

(b) Since in Theorem \ref{main} we proved in fact that
$W(\exp(-tS))\subseteq \Omega(\alpha)$, this means we give in the
present papers more precise localization of $W(\exp(-tS))$ than it
was claimed in \cite{CaZag}.

(c) The inclusion \eqref{CAHZAG} plays an important r\^{o}le for the
operator-norm error estimates of the norm convergence in the Euler
formula \eqref{Euler}, see the recent paper \cite{Zagrb1}, the
references quoted there and in particular \cite{Cach}. In the
present paper we give a new proof of the operator-norm convergent
Euler formula for the optimal error estimate with explicit
indication of its $\alpha-$dependence (\ref{alpha-estim}) for
$m$-$\alpha$-sectorial generators.

\vskip 1.5cm

\noindent \textbf{Acknowledgements}

\vskip 0.5cm

\noindent V.Z. is thankful to Pavel Bleher  for useful remarks
related to the multiplicative semigroups on the complex plane.


\end{document}